# CONSTRUCTION OF OPTIMAL MULTI-LEVEL SUPERSATURATED DESIGNS

By Hongquan Xu[1] and C. F. J. Wu[2]

*University of California, Los Angeles and Georgia Institute of Technology*

A supersaturated design is a design whose run size is not large enough for estimating all the main effects. The goodness of multi-level supersaturated designs can be judged by the generalized minimum aberration criterion proposed by Xu and Wu [*Ann. Statist.* **29** (2001) 1066–1077]. A new lower bound is derived and general construction methods are proposed for multi-level supersaturated designs. Inspired by the Addelman–Kempthorne construction of orthogonal arrays, several classes of optimal multi-level supersaturated designs are given in explicit form: Columns are labeled with linear or quadratic polynomials and rows are points over a finite field. Additive characters are used to study the properties of resulting designs. Some small optimal supersaturated designs of 3, 4 and 5 levels are listed with their properties.

**1. Introduction.** As science and technology have advanced to a higher level, investigators are becoming more interested in and capable of studying large-scale systems. Typically these systems have many factors that can be varied during design and operation. The cost of probing and studying a large-scale system can be prohibitively expensive. Building prototypes is time-consuming and costly. Even the quicker route of using computer modeling can take up many hours of CPU time. To address the challenges posed by this technological trend, research in experimental design has lately focused on the class of supersaturated designs for its run size economy and mathematical novelty. Formally, a supersaturated design (SSD) is a design whose run size is not large enough for estimating all the main effects represented by the columns of the design matrix. The design and analysis rely on

Received April 2003; revised January 2005.
[1]Supported in part by NSF Grant DMS-02-04009.
[2]Supported in part by NSF Grants DMS-00-72489 and DMS-03-05996.

the assumption of the effect sparsity principle ([5], [31], Section 3.5), that is, the number of relatively important effects in a factorial experiment is small. Some practical applications of SSDs can be found in [18, 19, 25, 30].

The construction of SSD dates back to Satterthwaite [26] and Booth and Cox [4]. The former suggested the use of random balanced designs and the latter proposed an algorithm to construct systematic SSDs. Many methods have been proposed for constructing two-level SSDs in the last decade, for example, among others, Lin [18, 19], Wu [30], Nguyen [25], Cheng [10], Li and Wu [16], Tang and Wu [29], Butler, Mead, Eskridge and Gilmour [7], Bulutoglu and Cheng [6] and Liu and Dean [20]. A popular criterion in the SSD literature is the $E(s^2)$ criterion [4], which measures the average correlation among columns. Nguyen [25] and Tang and Wu [29] independently derived the following lower bound for two-level SSDs with $N$ runs and $m$ factors:

$$(1) \qquad E(s^2) \geq N^2(m - N + 1)/[(m-1)(N-1)].$$

This lower bound was recently improved by Butler, Mead, Eskridge and Gilmour [7] and Bulutoglu and Cheng [6].

There are a few recent papers on the construction of multi-level SSDs. Yamada and Lin [37] and Yamada, Ikebe, Hashiguchi and Niki [36] proposed methods for the construction of three-level SSDs. Fang, Lin and Ma [13] and Lu and Sun [22] proposed algorithmic methods for the construction of multi-level SSDs. Lu, Hu and Zheng [21] constructed multi-level SSDs based on resolvable balanced incomplete block designs. Aggarwal and Gupta [2] proposed an algebraic construction method based on Galois field theory. Chatterjee and Gupta [8] studied multi-level SSDs with certain search design properties. Yamada and Matsui [38] considered the construction of mixed-level SSDs.

Extensions of the $E(s^2)$ criterion to the multi-level case are not unique, for example, ave$(\chi^2)$ statistic [37], ave$(f)$ [13] and $E(d^2)$ [22]. All these extensions measure the overall nonorthogonality between all possible pairs of columns. Lu and Sun [22] and Yamada and Matsui [38] derived lower bounds for $E(d^2)$ and ave$(\chi^2)$, respectively, which generalize the lower bound for $E(s^2)$ in (1).

There is another class of optimality criteria that were originally proposed for studying nonregular designs. *Generalized minimum aberration* (GMA) criteria, extensions of the *minimum aberration* criterion [14], have been proposed to assess regular or nonregular designs. See [12, 23, 28, 33, 34]. Obviously, these criteria can be used to assess SSDs as well. In this paper we adopt the GMA criterion due to Xu and Wu [34] as the optimality criterion for two reasons. First, the GMA criterion has good statistical properties and has been well justified for nonregular designs. Tang and Deng [28] and Xu and Wu [34] showed that GMA orthogonal designs are model robust in the sense that they tend to minimize the contamination of nonnegligible two-factor and higher-order interactions on the estimation of the main effects. Tang [27] and Ai and Zhang [3] provided projection justifications of the GMA criterion. Cheng, Deng and Tang [11] showed that GMA is also supported by some traditional model-dependent efficiency criteria. Second, the GMA criterion is general and can assess multi-level and mixed-level SSDs. It includes $E(s^2)$, ave($\chi^2$) and $E(d^2)$ as special cases. Section 2 reviews the GMA and other optimality criteria.

Section 3 presents some general optimality results for multi-level SSDs. We derive a new lower bound for multi-level SSDs as an extension of (1). This new lower bound is tight in many cases; for example, it is tight for SSDs with $s$ levels, $N = s^2$ runs and any number of factors. We also discuss the construction of optimal SSDs that achieve this bound. In particular, construction methods of Lin [18] and Tang and Wu [29] are extended to multi-level SSDs. Furthermore, optimal multi-level SSDs are shown to be periodic.

While an optimal SSD under GMA (and other criteria) may contain fully aliased columns, Section 4 describes construction methods that produce optimal SSDs without fully aliased columns. Inspired by the Addelman–Kempthorne construction of orthogonal arrays, we use linear and quadratic polynomials in the construction. Evaluating the polynomials over a finite field yields optimal multi-level SSDs. Compared with algorithmic methods, our algebraic method has at least two advantages: (i) the constructions are explicit and not limited to small run size, and (ii) the properties of resulting SSDs can be studied analytically. Section 5 presents some technical proofs that use additive characters of a finite field.

Section 6 lists some small optimal SSDs of 3, 4 and 5 levels and compares them with existing ones in terms of three other optimality criteria. For small run sizes ($\leq 25$), the benefits of our SSDs are marginal; our SSDs may be better in terms of one criterion but worse in terms of another criterion than existing ones. For 27 runs, our SSDs have more columns and one class of our designs is better than existing ones in terms of all three criteria. Section 7 considers mixed-level SSDs. A lower bound is derived and the method of replacement is proposed to construct optimal mixed-level SSDs.

**2. Optimality criteria.** Some definitions and notation are necessary in order to review the optimality criteria.

An $(N, s_1 s_2 \cdots s_m)$-design is an $N \times m$ matrix, where the $i$th column takes values on $s_i$ symbols $\{0, 1, \ldots, s_i - 1\}$. A design is called *mixed* if not all $s_i$'s are equal. Two designs are *isomorphic* if one can be obtained from the other through permutations of rows, columns and symbols in each column. An $OA(N, m, s, t)$ is an *orthogonal array* (OA) of $N$ runs, $m$ columns, $s$ levels and strength $t$, in which all possible level combinations appear equally often for any set of $t$ columns. An $OA(N, m, s, 1)$ is also called a balanced array.

We use the notation $SSD(N, s_1 s_2 \cdots s_m)$ to denote an SSD of $N$ runs, $m$ columns with levels $s_1, s_2, \ldots, s_m$ and use the notation $SSD(N, s^m)$ if all $s_i = s$. Throughout the paper, as in the literature, we only consider *balanced* SSDs, in which all levels appear equally often for any column.

2.1. *Generalized minimum aberration.* For an $(N, s_1 s_2 \cdots s_m)$-design $D$, consider the ANOVA model

$$Y = X_0 \alpha_0 + X_1 \alpha_1 + \cdots + X_m \alpha_m + \varepsilon,$$

where $Y$ is the vector of $N$ observations, $\alpha_j$ is the vector of all $j$-factor interactions, $X_j$ is the matrix of orthonormal contrast coefficients for $\alpha_j$, and $\varepsilon$ is the vector of independent random errors. For $j = 0, 1, \ldots, m$, Xu and Wu [34] defined $A_j(D)$, a function of $X_j$, to measure the aliasing between all $j$-factor interactions and the general mean. Specifically, if $X_j = [x_{ik}^{(j)}]$, define

$$A_j(D) = N^{-2} \sum_k \left| \sum_{i=1}^{N} x_{ik}^{(j)} \right|^2.$$

The GMA criterion sequentially minimizes the generalized wordlength patterns $A_1(D), A_2(D), A_3(D), \ldots$. Xu and Wu [34] showed that isomorphic designs have the same generalized wordlength patterns and therefore are not distinguishable under the GMA criterion.

The generalized wordlength patterns characterize the strength of a design, that is, $A_i(D) = 0$ for $i = 1, \ldots, t$ if and only if $D$ is an OA of strength $t$. For an SSD, $A_1 = 0$ and $A_2 > 0$. The GMA criterion suggests that we shall first minimize $A_2$ and then $A_3$, $A_4$ and so on. Note that $A_2$ measures the overall aliasing between all pairs of columns. Indeed, $A_2 = \sum_{i<j} r_{ij}^2$ if $R = (r_{ij})$ is the correlation matrix of all the main effects (see [32]). In particular, for a two-level design $A_2$ is equal to the sum of squares of correlations between all possible pairs of columns.

The following concept due to Xu [33] is useful in the theoretical development. For a design $D = [x_{ik}]$ and a positive integer $t$, define the $t$th power moment to be

$$K_t(D) = [N(N-1)/2]^{-1} \sum_{1 \le i < j \le N} [\delta_{ij}(D)]^t,$$

where $\delta_{ij}(D)$ is the number of coincidences between the $i$th and $j$th rows, that is, the number of $k$ such that $x_{ik} = x_{jk}$.

The following three lemmas are from [33].

LEMMA 1. *Suppose $D$ is an $SSD(N, s^m)$. Then*:

(i) $A_1(D) = 0$ and $A_2(D) = [(N-1)s^2 K_2(D) + m^2 s^2 - Nm(m+s-1)]/(2N)$;

(ii) $K_1(D) = m(N-s)/[(N-1)s]$ and $K_2(D) = [2NA_2(D) + Nm(m+s-1) - m^2 s^2]/[(N-1)s^2]$.

LEMMA 2. *Suppose $D$ is an $SSD(N, s^m)$. Then $A_2(D) \geq [m(s-1)(ms - m - N + 1)]/[2(N-1)]$.*

LEMMA 3. *Suppose $D$ is an $SSD(N, s^m)$. If the difference among all $\delta_{ij}(D)$, $i < j$, does not exceed one, then $D$ has GMA.*

2.2. *Other optimality criteria and connections.* There are several other optimality criteria for multi-level and mixed-level SSDs. Let $c_1, \ldots, c_m$ be the columns of an $SSD(N, s_1 s_2 \cdots s_m)$ and $n_{ab}^{ij}$ be the number of times that pair $(a, b)$ appears as a row in columns $c_i$ and $c_j$.

Yamada and Lin [37], Yamada, Ikebe, Hashiguchi and Niki [36] and Yamada and Matsui [38] defined

$$\chi^2(c_i, c_j) = \sum_{a=0}^{s_i-1} \sum_{b=0}^{s_j-1} [n_{ab}^{ij} - N/(s_i s_j)]^2 / (N/(s_i s_j))$$

to evaluate the dependency of columns $c_i$ and $c_j$. They proposed the following two criteria to evaluate the maximum and average dependency of all columns:

$$\mathrm{ave}(\chi^2) = \sum_{1 \leq i < j \leq m} \chi^2(c_i, c_j) / [m(m-1)/2]$$

and

$$\max(\chi^2) = \max_{1 \leq i < j \leq m} \chi^2(c_i, c_j).$$

Fang, Lin and Ma [13] defined

$$f(c_i, c_j) = \sum_{a=0}^{s_i-1} \sum_{b=0}^{s_j-1} |n_{ab}^{ij} - N/(s_i s_j)|$$

to measure the nonorthogonality between columns $c_i$ and $c_j$. They proposed to minimize

$$\mathrm{ave}(f) = \sum_{1 \leq i < j \leq m} f(c_i, c_j) / [m(m-1)/2] \quad \text{and} \quad \max(f) = \max_{1 \leq i < j \leq m} f(c_i, c_j)$$

and three other criteria.

When all $s_i = s$, Lu and Sun [22] and Lu, Hu and Zheng [21] defined

$$d_{ij}^2 = \sum_{a=0}^{s-1} \sum_{b=0}^{s-1} [n_{ab}^{ij} - N/s^2]^2$$

to measure the "departure from orthogonality" for columns $c_i$ and $c_j$. They proposed to minimize

$$E(d^2) = \sum_{1 \leq i < j \leq m} d_{ij}^2 / [m(m-1)/2] \quad \text{and} \quad \max(d^2) = \max_{1 \leq i < j \leq m} d_{ij}^2.$$

It is evident that $E(d^2) = (N/s^2)\operatorname{ave}(\chi^2)$ and $\max(d^2) = (N/s^2)\max(\chi^2)$.

LEMMA 4. *For an* $SSD(N, s^m)$:

(i) $\operatorname{ave}(\chi^2) = NA_2/[m(m-1)/2]$;
(ii) $E(d^2) = N^2 A_2/[s^2 m(m-1)/2]$.

The first part of Lemma 4 was proved by Xu [33] and the second part follows from the first part. Lemma 4 shows that $\operatorname{ave}(\chi^2)$ and $E(d^2)$ are equivalent to $A_2$; therefore, GMA can be viewed as a refinement of $\operatorname{ave}(\chi^2)$ and $E(d^2)$.

Note that $A_2$ measures the overall aliasing between columns. It is also important to measure the maximum aliasing between columns. For this purpose we consider projections and propose the concept of projected $A_2$ values. For a pair of columns $c_i$ and $c_j$, we define a *projected* $A_2$ value as $A_2(c_i, c_j) = A_2(d)$, where $d$ consists of the two columns $c_i$ and $c_j$. Obviously, the *overall* $A_2$ value is equal to the sum of all projected $A_2$ values, that is, $A_2(D) = \sum_{1 \leq i < j \leq m} A_2(c_i, c_j)$. Lemma 4 shows that the maximum projected $A_2$ value is equal to $\max(\chi^2)/N$ or $s^2 \max(d^2)/N^2$; therefore, the maximum projected $A_2$ value is equivalent to $\max(\chi^2)$ and $\max(d^2)$.

**3. Some optimality results.** In this and the next three sections we study multi-level SSDs.

3.1. *A new lower bound.*

THEOREM 1. *Suppose $D$ is an $SSD(N, s^m)$. Then*

$$A_2(D) \geq [m(s-1)(ms - m - N + 1)]/[2(N-1)] + (N-1)s^2 \eta(1-\eta)/(2N),$$

*where $\eta = m(N-s)/((N-1)s) - \lfloor m(N-s)/((N-1)s) \rfloor$ and $\lfloor x \rfloor$ is the largest integer that does not exceed $x$. The lower bound is achieved if and only if the number of coincidences, $\delta_{ij}(D)$, differs by at most one for all $i < j$. Furthermore, an SSD achieving the lower bound is optimal under GMA.*

PROOF. By Lemma 1(ii), $K_1(D) = [N(N-1)/2]^{-1} \sum_{1 \leq i < j \leq N} \delta_{ij}(D) = m(N-s)/((N-1)s)$. Then $\eta = K_1(D) - \lfloor K_1(D) \rfloor$ is the fractional part of $K_1(D)$. Since the number of coincidences, $\delta_{ij}(D)$, must be an integer, it is easy to verify that $K_2(D) = [N(N-1)/2]^{-1} \sum_{1 \leq i < j \leq N} [\delta_{ij}(D)]^2$ achieves the minimum value of $K_1(D)^2 + \eta(1-\eta)$ when all $\delta_{ij}(D)$ take on only one of the two values $\lfloor K_1(D) \rfloor$ and $\lfloor K_1(D) \rfloor + 1$. Then the lower bound of $A_2(D)$ follows from Lemma 1(i) and some straightforward algebra. By Lemma 3 an SSD achieving this lower bound is optimal under GMA. □

The lower bound in Theorem 1 is new for multi-level SSDs. Using the connection established in Lemma 4, we obtain a new lower bound for $\mathrm{ave}(\chi^2)$ and $E(d^2)$. As will be seen next, there are many cases in which the lower bound in Theorem 1 is tight, whereas the lower bound in Lemma 2 is not. For example, when $N = s^2$, the lower bound in Theorem 1 is tight for any $m$; in contrast, the lower bound in Lemma 2 is tight only when $m$ is a multiple of $s+1$.

3.2. *Optimal designs.* Many optimal SSDs that achieve the lower bound in Theorem 1 can be derived from saturated OAs. An $OA(N, t, s, 2)$ is saturated if $N - 1 = t(s-1)$. The following lemma from Mukerjee and Wu [24] says that the number of coincidences between distinct rows is constant for a saturated OA.

LEMMA 5. *Suppose $H$ is a saturated $OA(N, t, s, 2)$ with $t = (N-1)/(s-1)$. Then $\delta_{ij}(H) = (N-s)/[s(s-1)]$ for any $i < j$.*

The next lemma shows the change of the $A_2$ values of a design when a saturated OA is juxtaposed to it. Readers are referred to Xu and Wu [35] for a proof.

LEMMA 6. *Suppose $H$ is a saturated $OA(N, t, s, 2)$ with $t = (N-1)/(s-1)$ and $D$ is an $OA(N, m, s, 1)$. Let $D \cup H$ be the column juxtaposition of $D$ and $H$. Then $A_2(D \cup H) = A_2(D) + m(s-1)$.*

Tang and Wu [29] first proposed to construct optimal two-level SSDs by juxtaposing saturated OAs derived from Hadamard matrices. This method can be extended to construct optimal multi-level SSDs. Suppose $D_1, \ldots, D_k$ are $k$ saturated $OA(N, t, s, 2)$s with $t = (N-1)/(s-1)$. Let $D = D_1 \cup \cdots \cup D_k$ be the *column* juxtaposition, which may have duplicated or fully aliased columns. It is evident that $\delta_{ij}(D) = k(N-s)/[s(s-1)]$ for any $i < j$. Then by Lemma 3, $D$ is an optimal SSD under GMA.

As Tang and Wu [29] suggested, to construct an SSD with $m = kt - j$ columns, $1 \le j < t$, we may simply delete the last $j$ columns from $D$. Though the resulting design may not be optimal, it has an $A_2$ value very close to the lower bound in Theorem 1.

If one column is removed from or one (balanced) column is added to $D$, the resulting design is still optimal. Cheng [10] showed that, for two-level

SSDs, removing (and resp. adding) two orthogonal columns from (and resp. to) $D$ also results in an optimal SSD. This is not true for multi-level SSDs in general. For $N = s^2$, we have a stronger result in Lemma 5 that the number of coincidences between any two rows is equal to 1. Then removing (and resp. adding) any number of orthogonal columns from (and resp. to) $D$ also results in an optimal SSD under GMA, because the resulting design has the property that the number of coincidences between any two rows differs by at most one. In particular, for any $m$, the lower bound in Theorem 1 is tight.

Lin [18] used half fractions of Hadamard matrices to construct two-level SSDs by taking a column as the branching column. This method can be extended to construct multi-level SSDs as follows. Taking any column of a saturated $OA(N, t, s, 2)$ as the branching column, we obtain $s$ fractions according to the levels of the branching column. After removing the branching column, the fractions have the properties that all columns are balanced and the number of coincidences between any two rows is constant. The *row juxtaposition* of any $k$ fractions forms an $SSD(kNs^{-1}, s^{t-1})$ of which the number of coincidences between any two rows differs by at most one. By Lemma 3 such a design is optimal under GMA. For $N = s^2$, any subdesign is also optimal, because the number of coincidences between any two rows is either 0 or 1.

Because a saturated $OA(s^n, (s^n - 1)/(s - 1), s, 2)$ exists for any prime power $s$, we have the following result.

THEOREM 2. *Suppose $s$ is a prime power.*

(i) *For any $n$ and $k$, there exists an optimal $SSD(s^n, s^m)$ that achieves the lower bound in Theorem 1, where $m = k(s^n - 1)/(s - 1)$ or $m = k(s^n - 1)/(s - 1) \pm 1$.*

(ii) *For any $n$ and $k < s$, there exists an optimal $SSD(ks^{n-1}, s^m)$ that achieves the lower bound in Theorem 1, where $m = (s^n - 1)/(s - 1) - 1$.*

(iii) *For any $m$, there exists an optimal $SSD(s^2, s^m)$ that achieves the lower bound in Theorem 1.*

(iv) *For any $m \leq s$ and $k < s$, there exists an optimal $SSD(ks, s^m)$ that achieves the lower bound in Theorem 1.*

Given $N$ and $s$, let $a_2(m) = \min\{A_2(D) : D \text{ is an } SSD(N, s^m)\}$, where designs may have fully aliased columns. When $N = s^2$, Theorem 2(iii) implies that $a_2(m + s + 1) = a_2(m) + m(s - 1)$ for any $m \geq 1$. The following result shows that for certain $N$, $a_2(m)$ is periodic when $m$ is large enough.

THEOREM 3. *Suppose a saturated $OA(N, t, s, 2)$ exists with $t = (N - 1)/(s - 1)$. Then there exists a positive integer $m_0$ such that for $m \geq m_0$, $a_2(m + t) = a_2(m) + m(s - 1)$.*

Readers are referred to Xu and Wu [35] for a proof. This periodicity property helps us understand SSDs of large size; it shows how larger optimal SSDs are connected with smaller ones. Chen and Wu [9] previously showed a similar periodicity property for maximum resolution and minimum aberration designs.

The above optimal SSDs may contain fully aliased columns, which are not useful in practice. The next section presents explicit construction methods that produce optimal SSDs without fully aliased columns.

**4. Construction.** The construction methods are applicable to any prime power. Throughout this section we assume $s > 2$ is a prime power. Let $F_s$ be a Galois field of $s$ elements. For clarity, all proofs are given in the next section.

4.1. *Half Addelman–Kempthorne orthogonal arrays.* Addelman and Kempthorne [1] described a method for constructing $OA(2s^n, 2(s^n - 1)/(s - 1) - 1, s, 2)$ for any prime power $s$ and any $n$. Such arrays can be naturally decomposed into two arrays of $s^n$ runs. Each array is an $SSD(s^n, s^m)$ with $m = 2(s^n - 1)/(s - 1) - 1$. We now describe how to construct an SSD in general.

In the construction the columns of an array are labeled with linear or quadratic polynomials in $n$ variables $X_1, \ldots, X_n$ and the rows are labeled with points from $F_s^n$. Let $f_1(X_1, \ldots, X_n)$ and $f_2(X_1, \ldots, X_n)$ be two functions, linear or nonlinear. They correspond to two columns of length $s^n$ when evaluated at $F_s^n$. The two functions (or columns) are *fully aliased* if the pair has only $s$ level combinations, each combination occurring $s^{n-1}$ times; and *orthogonal* if the pair has $s^2$ distinct level combinations, each combination occurring $s^{n-2}$ times. A pair of fully aliased columns has projected $A_2 = s - 1$ and a pair of orthogonal columns has projected $A_2 = 0$.

Following Addelman and Kempthorne [1], $f_1(X_1, \ldots, X_n)$ and $f_2(X_1, \ldots, X_n)$ are said to be *semi-orthogonal* to each other if (i) for $s$ odd, the pair has $(s+1)s/2$ distinct level combinations, $s$ combinations occurring $s^{n-2}$ times and $s(s-1)/2$ combinations occurring $2s^{n-2}$ times, and (ii) for $s$ even, the pair has $s^2/2$ distinct level combinations each occurring $2s^{n-2}$ times. A pair of semi-orthogonal columns has projected $A_2 = (s-1)/s$ for $s$ odd and projected $A_2 = 1$ for $s$ even. This result can be easily verified from the connection between the ave($\chi^2$) statistic and $A_2$ described in Lemma 4.

Let $L(X_1, \ldots, X_n)$ be the set of all nonzero linear functions of $X_1, \ldots, X_n$, that is,

$$L(X_1, \ldots, X_n) = \{c_1 X_1 + \cdots + c_n X_n : c_i \in F_s, \text{ not all } c_i \text{ are zero}\}.$$

Every function in $L(X_1, \ldots, X_n)$ corresponds to a balanced column. Two functions $f_1$ and $f_2$ in $L(X_1, \ldots, X_n)$ are *dependent* if there is a nonzero constant $c \in F_s$ such that $f_1 = cf_2$; otherwise, they are *independent*. Clearly, dependent linear functions correspond to the same column up to level permutation and, thus, they are fully aliased, while independent linear functions correspond to orthogonal columns. A set of $(s^n - 1)/(s - 1)$ independent linear functions generates an $OA(s^n, (s^n - 1)/(s - 1), s, 2)$. The traditional convention is to assume the first nonzero element is 1 for each column. For convenience, we assume the *last* nonzero element is 1 for each column. In particular, let $H(X_1, \ldots, X_n)$ be the set of all nonzero linear functions of $X_1, \ldots, X_n$ such that the *last* nonzero coefficient is 1. When evaluated at $F_s^n$, $H(X_1, \ldots, X_n)$ is a saturated $OA(s^n, (s^n - 1)/(s - 1), s, 2)$. This is indeed the regular fractional factorial design and the construction is called the Rao–Hamming construction by Hedayat, Sloane and Stufken ([15], Section 3.4).

The key idea of the Addelman–Kempthorne construction is to use quadratic functions in addition to linear functions. Let

(2) $\quad Q_1^*(X_1, \ldots, X_n) = \{X_1^2 + aX_1 + h : a \in F_s, h \in H(X_2, \ldots, X_n)\}$

and $Q_1(X_1, \ldots, X_n) = \{X_1\} \cup Q_1^*(X_1, \ldots, X_n)$.

$H(X_1, \ldots, X_n)$ has $(s^n - 1)/(s - 1)$ columns and $Q_1^*(X_1, \ldots, X_n)$ has $(s^n - 1)/(s - 1) - 1$ columns. The column juxtaposition of $H(X_1, \ldots, X_n)$ and $Q_1^*(X_1, \ldots, X_n)$ forms an $SSD(s^n, s^m)$ with $m = 2(s^n - 1)/(s - 1) - 1$, which is a half of an Addelman–Kempthorne OA.

EXAMPLE 1. Consider $s = 3$ and $n = 2$. The functions are

$$H(X_1, X_2) = \{X_1, X_2, X_1 + X_2, 2X_1 + X_2\},$$
$$Q_1^*(X_1, X_2) = \{X_1^2 + X_2, X_1^2 + X_1 + X_2, X_1^2 + 2X_1 + X_2\},$$
$$Q_1(X_1, X_2) = \{X_1, X_1^2 + X_2, X_1^2 + X_1 + X_2, X_1^2 + 2X_1 + X_2\}.$$

$H(X_1, X_2)$ is an $OA(9, 4, 3, 2)$ when the functions are evaluated at $F_3^2$; so is $Q_1(X_1, X_2)$. They are isomorphic [indeed, there is only one unique $OA(9, 4, 3, 2)$ up to isomorphism]. The column juxtaposition of $H(X_1, X_2)$ and $Q_1^*(X_1, X_2)$ forms an $SSD(9, 3^7)$, which is isomorphic to the first (and last) nine rows of the commonly used $OA(18, 7, 3, 2)$ (e.g., Table 7C.2 of [31]). This SSD has an overall $A_2 = 6$ and achieves the lower bound in Theorem 1. Furthermore, there are no fully aliased columns. Each column of $Q_1^*(X_1, X_2)$ is semi-orthogonal to three columns of $H(X_1, X_2)$ with projected $A_2 = 2/3$.

In general, we have the following results.

LEMMA 7. *When evaluated at $F_s^n$, $Q_1(X_1, \ldots, X_n)$ is an $OA(s^n, (s^n - 1)/(s - 1), s, 2)$.*

THEOREM 4. *The column juxtaposition of $H(X_1, \ldots, X_n)$ and $Q_1^*(X_1, \ldots, X_n)$ forms an $SSD(s^n, s^m)$ with $m = 2(s^n - 1)/(s - 1) - 1$. It has an overall $A_2 = s^n - s$ and is optimal under GMA. Furthermore, column $X_1$ is orthogonal to all other columns and there are no fully aliased columns if $s > 2$.*

(i) *For $s$ odd, the possible projected $A_2$ values are $0$ and $(s-1)/s$. There are $s(s^n - s)/(s - 1)$ pairs of semi-orthogonal columns with projected $A_2 = (s-1)/s$.*

(ii) *For $s$ even, the possible projected $A_2$ values are $0$ and $1$. There are $s^n - s$ pairs of semi-orthogonal columns with projected $A_2 = 1$.*

Both $Q_1(X_1, \ldots, X_n)$ and $H(X_1, \ldots, X_n)$ are saturated OAs of the same parameters. It is of interest to know whether they are isomorphic. Example 1 shows that they are isomorphic for $n = 2$ and $s = 3$. This is true as long as $n = 2$. When $n > 2$ and $s > 2$, they are not isomorphic. The following corollary summarizes the result.

COROLLARY 1. (i) *For $n = 2$, $Q_1(X_1, X_2)$ is isomorphic to the regular design $H(X_1, X_2)$.*

(ii) *For $n > 2$ and $s > 2$, $Q_1(X_1, \ldots, X_n)$ is not isomorphic to $H(X_1, \ldots, X_n)$.*

Corollary 1(ii) implies that $Q_1(X_1, \ldots, X_n)$ is a nonregular design for $n > 2$ and $s > 2$.

4.2. *Juxtaposition of saturated orthogonal arrays.* As a by-product of the half Addelman–Kempthorne construction, we have constructed a saturated OA, $Q_1(X_1, \ldots, X_n)$, besides the regular OA, $H(X_1, \ldots, X_n)$. For any $h \in H(X_1, \ldots, X_n)$, we can construct a saturated OA, $Q_h(X_1, \ldots, X_n)$, as follows. Let $h = c_1 X_1 + \cdots + c_n X_n$, not all $c_i$'s being $0$. Let $k$ be the largest $i$ such that $c_i \neq 0$. Then $c_k = 1$ and $c_i = 0$ for all $i > k$. Let $Y_1 = h$, $Y_i = X_{i-1}$ for $2 \leq i \leq k$, and $Y_i = X_i$ for $k < i \leq n$. It is clear that $H(X_1, \ldots, X_n)$ is equivalent to $H(Y_1, \ldots, Y_n)$ up to row and column permutations. Define $Q_h^*(X_1, \ldots, X_n) = Q_1^*(Y_1, \ldots, Y_n)$ as in (2) by replacing $X_i$ with $Y_i$ and $Q_h(X_1, \ldots, X_n) = Q_1(Y_1, \ldots, Y_n)$.

Since there are $(s^n - 1)/(s - 1)$ columns in $H(X_1, \ldots, X_n)$, we obtain $(s^n - 1)/(s - 1)$ saturated $OA(s^n, (s^n - 1)/(s - 1), s, 2)$s. Although they are all isomorphic, we can obtain many optimal multi-level SSDs by juxtaposing them.

EXAMPLE 2. Consider $s = 3$ and $n = 2$. $H(X_1, X_2) = \{X_1, X_2, X_1 + X_2, 2X_1 + X_2\}$. For each $h \in H(X_1, X_2)$, we can define $Q_h(X_1, X_2)$ as follows:
$$Q_{X_1}(X_1, X_2) = \{X_1, X_1^2 + X_2, X_1^2 + X_1 + X_2, X_1^2 + 2X_1 + X_2\},$$

$$Q_{X_2}(X_1, X_2) = \{X_2, X_2^2 + X_1, X_2^2 + X_2 + X_1, X_2^2 + 2X_2 + X_1\},$$
$$Q_{X_1+X_2}(X_1, X_2) = \{X_1 + X_2, (X_1 + X_2)^2 + X_1,$$
$$(X_1 + X_2)^2 + 2X_1 + X_2, (X_1 + X_2)^2 + 2X_2\},$$
$$Q_{2X_1+X_2}(X_1, X_2) = \{2X_1 + X_2, (2X_1 + X_2)^2 + X_1,$$
$$(2X_1 + X_2)^2 + X_2, (2X_1 + X_2)^2 + 2X_1 + 2X_2\}.$$

Each $Q_h(X_1, X_2)$ is a saturated $OA(9, 4, 3, 2)$ and they are all isomorphic. The column juxtaposition of all four $Q_h(X_1, X_2)$ has 16 columns: 4 linear and 12 quadratic. All linear columns are orthogonal to each other. Each linear column is orthogonal to 3 quadratic columns, and semi-orthogonal to the other 9 quadratic columns. Each quadratic column is orthogonal to 1 linear column, semi-orthogonal to the other 3 linear columns, orthogonal to 2 quadratic columns, and partially aliased (projected $A_2 = 4/9$) with the other 9 quadratic columns. The 16 columns together form an optimal $SSD(9, 3^{16})$ with an overall $A_2 = 48$. The 12 quadratic columns together form an optimal $SSD(9, 3^{12})$ with an overall $A_2 = 24$. For the latter design, each column is partially aliased with 9 columns with projected $A_2 = 4/9$.

THEOREM 5. *Let $h_1$ and $h_2$ be two distinct functions in $H(X_1, \ldots, X_n)$. The column juxtaposition of $Q_{h_1}(X_1, \ldots, X_n)$ and $Q_{h_2}(X_1, \ldots, X_n)$ forms an $SSD(s^n, s^m)$ with $m = 2(s^n - 1)/(s - 1)$. It has an overall $A_2 = s^n - 1$ and is optimal under GMA. Furthermore, there are no fully aliased columns if $s$ is odd or $s > 4$.*

*(i) For $s$ odd, the possible projected $A_2$ values are $0$, $(s-1)/s$, $(s-1)^2/s^2$ and $(s-1)/s^2$. There are $2s$ pairs with projected $A_2 = (s-1)/s$, $s^2$ pairs with projected $A_2 = (s-1)^2/s^2$ and $s^2(s^n - s^2)/(s-1)$ pairs with projected $A_2 = (s-1)/s^2$.*

*(ii) For $s$ even, the possible projected $A_2$ values are $0, 1, 2$ and $3$.*

*(iii) For $s = 4$, the possible projected $A_2$ values are $0, 1$ and $3$. There are one pair of fully aliased columns with projected $A_2 = 3$ and $4^n - 4$ pairs of partially aliased columns with projected $A_2 = 1$.*

Theorem 5 states that the column juxtaposition of $Q_{h_1}(X_1, \ldots, X_n)$ and $Q_{h_2}(X_1, \ldots, X_n)$ has the same projected $A_2$ values and frequencies, independent of the choice of $h_1$ and $h_2$. It is of interest to note that they can have different geometric structures and be nonisomorphic to each other. For example, when $n = 3$ and $s = 3$, the column juxtaposition of $Q_{X_1}$ and $Q_{X_2}$ is not isomorphic to the column juxtaposition of $Q_{X_1}$ and $Q_{X_3}$.

Extending Theorem 5, we have the following result.

THEOREM 6. *For $1 < k \leq (s^n - 1)/(s - 1)$, let $h_1, \ldots, h_k$ be $k$ distinct functions in $H(X_1, \ldots, X_n)$. The column juxtaposition of $Q_{h_i}(X_1, \ldots, X_n)$, $i = 1, \ldots, k$, forms an $SSD(s^n, s^m)$ with $m = k(s^n - 1)/(s - 1)$. It has an overall $A_2 = \binom{k}{2}(s^n - 1)$ and is optimal under GMA. Furthermore, there are no fully aliased columns if $s$ is odd or $s > 4$.*

*(i) For $s$ odd, the possible projected $A_2$ values are $0$, $(s-1)/s$, $(s-1)^2/s^2$ and $(s-1)/s^2$. There are $\binom{k}{2}2s$ pairs with projected $A_2 = (s-1)/s$, $\binom{k}{2}s^2$ pairs with projected $A_2 = (s-1)^2/s^2$ and $\binom{k}{2}s^2(s^n - s^2)/(s-1)$ pairs with projected $A_2 = (s-1)/s^2$.*

*(ii) For $s$ even, the possible projected $A_2$ values are $0, 1, 2$ and $3$.*

*(iii) For $s = 4$, the possible projected $A_2$ values are $0, 1$ and $3$. There are $\binom{k}{2}$ pairs of fully aliased columns with projected $A_2 = 3$ and $\binom{k}{2}(4^n - 4)$ pairs of partially aliased columns with projected $A_2 = 1$.*

When $k = (s^n - 1)/(s-1)$, the above SSD has $[(s^n - 1)/(s-1)]^2$ columns, among which $(s^n - 1)/(s - 1)$ columns are linear from $H(X_1, \ldots, X_n)$ and the rest are quadratic. All quadratic functions form another class of SSDs. This SSD does not have semi-orthogonal columns, which have projected $A_2 = (s-1)/s$ for $s$ odd.

THEOREM 7. *Suppose $s$ is odd. For $1 < k \leq (s^n - 1)/(s - 1)$, let $h_1, \ldots, h_k$ be $k$ distinct functions in $H(X_1, \ldots, X_n)$. The column juxtaposition of $Q_{h_i}^*(X_1, \ldots, X_n)$, $i = 1, \ldots, k$, forms an $SSD(s^n, s^m)$ with $m = k(s^n - s)/(s - 1)$. There are no fully aliased columns and the possible projected $A_2$ values are $0$, $(s-1)^2/s^2$ and $(s-1)/s^2$. There are $\binom{k}{2}s^2$ pairs with projected $A_2 = (s-1)^2/s^2$ and $\binom{k}{2}s^2(s^n - s^2)/(s-1)$ pairs with projected $A_2 = (s-1)/s^2$. It has an overall $A_2 = \binom{k}{2}(s^n - 2s + 1)$. When $k = (s^n - 1)/(s-1) - 1$ or $(s^n - 1)/(s - 1)$, the SSD is optimal under GMA.*

COROLLARY 2. *For $s$ odd, the column juxtaposition of $Q_h^*(X_1, X_2)$, $h \in H(X_1, X_2)$, forms an $SSD(s^2, s^{(s+1)s})$. It has an overall $A_2 = (s+1)s(s-1)^2/2$ and is optimal under GMA. Each column is orthogonal to $s-1$ columns and partially aliased with the other $s^2$ columns with projected $A_2 = (s-1)^2/s^2$.*

4.3. *Fractions of saturated orthogonal arrays.* First consider fractions of $H(X_1, \ldots, X_n)$. Without loss of generality, taking $X_1$ as the branching column, we obtain $s$ fractions according to the levels of $X_1$. Each fraction has $s^{n-1}$ runs and $(s^n - 1)/(s - 1)$ columns: column $X_1$ has only one level and all other columns have $s$ levels. The row juxtaposition of any $k$ fractions forms an optimal SSD after removing column $X_1$.

THEOREM 8. *Take any column of $H(X_1,\ldots,X_n)$ as a branching column. For $k < s$, the row juxtaposition of any $k$ fractions forms an $SSD(ks^{n-1}, s^m)$ with $m = (s^n - s)/(s - 1)$ after removing the branching column. It has an overall $A_2 = (s^n - s)(s - k)/(2k)$ and is optimal under GMA. Furthermore, all possible projected $A_2$ values are $0$ and $(s - k)/k$. There are $(s^n - s)/2$ pairs of nonorthogonal columns with projected $A_2 = (s - k)/k$. In particular, there are no fully aliased columns for $1 < k < s$.*

Next consider fractions of $Q_1(X_1,\ldots,X_n)$. If $X_1$ is used as the branching column, the row juxtaposition of the fractions has the same property as that of $H(X_1,\ldots,X_n)$. In the following theorem, we take $X_1^2 + X_2$ as the branching column.

THEOREM 9. *Take column $X_1^2 + X_2$ of $Q_1(X_1,\ldots,X_n)$ as a branching column. The row juxtaposition of any $k$ fractions forms an $SSD(ks^{n-1}, s^m)$ with $m = (s^n - s)/(s - 1)$ after removing the branching column. It has an overall $A_2 = (s^n - s)(s - k)/(2k)$ and is optimal under GMA. Furthermore, there are no fully aliased columns for $1 < k < s$.*

*(i) For $s$ odd, there are $s(s^n - s^2 + s - 1)/2$ pairs of nonorthogonal columns, $s(s - 1)/2$ pairs with projected $A_2 = (s - k)/k$ and $s(s^n - s^2)/2$ pairs with projected $A_2 = (s - k)/(ks)$.*

*(ii) For $s$ even, there are at most $(s-1)(s^n - s^2 + s)/2$ pairs of nonorthogonal columns, $s(s - 1)/2$ pairs with projected $A_2 = (s - k)/k$ and at most $(s - 1)(s^n - s^2)/2$ pairs with projected $A_2 \leq 1$.*

*(iii) For $s = 4$ and $k = 2$, there are $(4^n - 4)/2$ pairs of nonorthogonal columns with projected $A_2 = 1$; for $s = 4$ and $k = 3$, there are $6$ pairs of nonorthogonal columns with projected $A_2 = 1/3$ and $3(4^n - 16)/2$ pairs with projected $A_2 = 1/9$.*

By branching other columns, we can obtain different SSDs as illustrated below.

EXAMPLE 3. Consider $n = 3$ and $s = 3$. The columns of $Q_1(X_1, X_2, X_3)$ are the following:

$$X_1, X_1^2 + X_2, X_1^2 + X_1 + X_2, X_1^2 + 2X_1 + X_2, X_1^2 + X_3, X_1^2 + X_1 + X_3,$$

$$X_1^2 + 2X_1 + X_3, X_1^2 + X_2 + X_3, X_1^2 + X_1 + X_2 + X_3, X_1^2 + 2X_1 + X_2 + X_3,$$

$$X_1^2 + 2X_2 + X_3, X_1^2 + X_1 + 2X_2 + X_3, X_1^2 + 2X_1 + 2X_2 + X_3.$$

Depending on the branching column, we obtain one of three types of optimal $SSD(18,3^{12})$. The frequencies of projected $A_2$ values are the following:

| $A_2$ | 0 | 1/6 | 1/2 |
|---|---|---|---|
| Type 1 | 54 | 0 | 12 |
| Type 2 | 36 | 27 | 3 |
| Type 3 | 42 | 18 | 6 |

We obtain a type 1 SSD if $X_1$ is used as the branching column, a type 2 SSD if $X_1^2 + aX_1 + X_2$ is used as the branching column, and a type 3 SSD if $X_1^2 + aX_1 + bX_2 + X_3$ is used as the branching column, where $a, b \in F_3$. A type 2 design is preferred in general because it has the smallest number of maximum projected $A_2$.

**5. Some proofs.** Additional notation and lemmas are needed for the proofs. Let $\mathcal{C}$ be the set of complex numbers and $F_s^*$ be the set of nonzero elements in $F_s$. An additive *character* of $F_s$ is a homomorphism mapping $\chi: F_s \to \mathcal{C}$ such that, for any $x, y \in F_s$, $|\chi(x)| = 1$ and $\chi(x+y) = \chi(x)\chi(y)$. Clearly, $\chi(0) = 1$ since $\chi(0) = \chi(0)\chi(0)$. A character is called *trivial* if $\chi(x) = 1$ for all $x$; otherwise, it is *nontrivial*. A nontrivial additive character has the property that, for $a \in F_s$, $\sum_{x \in F_s} \chi(ax) = s$ if $a = 0$ and equals 0 otherwise (see, e.g., [17]).

Let $\chi$ be a nontrivial additive character. For $u \in F_s$, the function $\chi_u(x) = \chi(ux)$ defines a character of $F_s$. Then $\chi_0$ is a trivial character and all other characters $\chi_u$ are nontrivial. It is important to note that $\{\chi_u, u \in F_s^*\}$ forms a set of orthonormal contrasts defined in [34], that is, $\sum_{x \in F_s} \chi_u(x)\overline{\chi_v(x)} = s$ if $u = v$ and equals 0 otherwise. As a result, we can use additive characters to compute the generalized wordlength pattern. In particular, for a column $x = (x_1, \ldots, x_N)^T$, the orthonormal contrast coefficient matrix is $(\chi_u(x_i))$, where $u \in F_s^*$ and $i = 1, \ldots, N$. For a pair of columns $x = (x_1, \ldots, x_N)^T$ and $y = (y_1, \ldots, y_N)^T$, the projected $A_2$ value is

$$(3) \qquad A_2(x,y) = N^{-2} \sum_{u_1 \in F_s^*} \sum_{u_2 \in F_s^*} \left| \sum_{i=1}^{N} \chi(u_1 x_i + u_2 y_i) \right|^2.$$

Let $s = p^r$, where $p$ is a prime. Define a mapping $\mathrm{Tr}: F_s \to F_p$, called the *trace*, as follows: $\mathrm{Tr}(x) = x + x^p + x^{p^2} + \cdots + x^{p^{r-1}}$ for any $x \in F_s$. Let

$$(4) \qquad \chi(x) = e^{2\pi i \, \mathrm{Tr}(x)/p} \qquad \text{for any } x \in F_s.$$

This is a nontrivial additive character and is called the *canonical* additive character of $F_s$. For $s = 2$, the canonical additive character is the usual contrast coding: $\chi(0) = 1$ and $\chi(1) = -1$.

The following three lemmas are useful when evaluating projected $A_2$ values. Interested readers are referred to Xu and Wu [35] for proofs of the lemmas in this section.

LEMMA 8. *For $s$ odd, let $a \in F_s^*$, $b, c \in F_s$, and $\chi$ be a nontrivial additive character. Then $|\sum_{x \in F_s} \chi(ax^2 + bx + c)|^2 = s$.*

LEMMA 9. *For $s$ even, let $a, b \in F_s$ and $\chi$ be the canonical additive character of $F_s$ defined in (4). Then $\sum_{x \in F_s} \chi(ax^2 + bx) = s$ if $a = b^2$ and equals 0 otherwise.*

LEMMA 10. *Let $G$ be a subset of $F_s$, $|G| = k$ and $\chi$ be a nontrivial additive character. Then $\sum_{u \in F_s^*} |\sum_{x \in G} \chi(ux)|^2 = (s-k)k$.*

The following lemma follows from Lemmas 1–4 and 5a of Addelman and Kempthorne [1].

LEMMA 11. *Consider columns $X_1^2 + a_1 X_1 + h_1$ and $a_2 X_1 + h_2$, where $h_1, h_2 \in L(X_2, \ldots, X_n)$ and $a_1, a_2 \in F_s$.*

(i) *If $h_1$ and $h_2$ are independent, they are orthogonal.*
(ii) *For $s$ odd, if $h_1$ and $h_2$ are dependent, they are semi-orthogonal.*
(iii) *For $s$ even, if $h_1$ and $h_2$ are dependent and $a_1 h_2 = a_2 h_1$, they are orthogonal.*
(iv) *For $s$ even, if $h_1$ and $h_2$ are dependent and $a_1 h_2 \neq a_2 h_1$, they are semi-orthogonal.*

PROOF OF THEOREM 4. The columns of $Q_1(X_1, \ldots, X_n)$ are $X_1$ and $X_1^2 + a_1 X_1 + h_1$, and the columns of $H(X_1, \ldots, X_n)$ are $X_1$ and $a_2 X_1 + h_2$, where $a_i \in F_s$ and $h_i \in H(X_2, \ldots, X_n)$. Since both $H(X_1, \ldots, X_n)$ and $Q_1(X_1, \ldots, X_n)$ are saturated OAs and they share column $X_1$, the optimality of the column juxtaposition of $H(X_1, \ldots, X_n)$ and $Q_1^*(X_1, \ldots, X_n)$ follows from Lemmas 3 and 5. By Lemma 6, the overall $A_2(H \cup Q_1^*) = A_2(Q_1^*) + [(s^n - s)/(s-1)](s-1) = (s^n - s)$ since $Q_1^*$ is an OA of strength 2.

(i) When $s$ is odd, by Lemma 11, $X_1^2 + a_1 X_1 + h_1$ and $a_2 X_1 + h_2$ are semi-orthogonal if $h_1 = h_2$. Therefore, each column of $Q_1^*(X_1, \ldots, X_n)$ is semi-orthogonal to $s$ columns of $H(X_1, \ldots, X_n)$. Since there are $(s^n - 1)/(s-1) - 1$ columns in $Q_1^*(X_1, \ldots, X_n)$, there are in total $s(s^n - s)/(s-1)$ semi-orthogonal pairs of columns with projected $A_2 = (s-1)/s$.

(ii) When $s$ is even, by Lemma 11, $X_1^2 + a_1 X_1 + h_1$ and $a_2 X_1 + h_2$ are semi-orthogonal if $h_1 = h_2$ and $a_1 \neq a_2$. Therefore, each column of $Q_1^*(X_1, \ldots, X_n)$ is semi-orthogonal to $s - 1$ columns of $H(X_1, \ldots, X_n)$. Since there are $(s^n - 1)/(s-1) - 1$ columns in $Q_1^*(X_1, \ldots, X_n)$, there are in total $s^n - s$ semi-orthogonal pairs of columns with projected $A_2 = 1$. □

PROOF OF COROLLARY 1. (i) Let $Y_1 = X_1$ and $Y_2 = X_1^2 + X_2$. This is a one-to-one mapping from $(Y_1, Y_2)$ to $(X_1, X_2)$. The columns of $Q_1(X_1, X_2)$ are $X_1 = Y_1$ and $X_1^2 + aX_1 + X_2 = aY_1 + Y_2$, where $a \in F_s$. Therefore, $Q_1(X_1, X_2) = H(Y_1, Y_2)$ is isomorphic to $H(X_1, X_2)$.

(ii) It follows from Theorems 8 and 9 to be proven later. □

LEMMA 12. *Suppose $h_i \in L(X_3, \ldots, X_n)$ and $a_i, b_i \in F_s$ for $i = 1, 2$.*

(i) *If $h_1$ and $h_2$ are independent, $X_1^2 + a_1X_1 + b_1X_2 + h_1$ and $X_2^2 + a_2X_2 + b_2X_1 + h_2$ are orthogonal.*

(ii) *If $b_2 \neq 0$, $X_1^2 + a_1X_1 + b_1X_2 + h_1$ and $X_2^2 + a_2X_2 + b_2X_1$ are orthogonal.*

(iii) *If $h_1$ and $h_2$ are dependent, the pair of columns $X_1^2 + a_1X_1 + b_1X_2 + h_1$ and $X_2^2 + a_2X_2 + b_2X_1 + h_2$ has projected $A_2 = (s-1)/s^2$ for $s$ odd and $A_2 = 0$ or $1$ for $s$ even.*

(iv) *For $s$ odd, the pair of columns $X_1^2 + a_1X_1 + X_2$ and $X_2^2 + a_2X_2 + X_1$ has projected $A_2 = (s-1)^2/s^2$.*

(v) *For $s$ even, the pair of columns $X_1^2 + a_1X_1 + X_2$ and $X_2^2 + a_2X_2 + X_1$ has projected $A_2 = 0, 1, 2$ or $3$.*

(vi) *For $s = 4$, the pair of columns $X_1^2 + a_1X_1 + X_2$ and $X_2^2 + a_2X_2 + X_1$ has projected $A_2 = 3$ if $a_1 = a_2 = 0$, $A_2 = 1$ if both $a_1 \neq 0$ and $a_2 \neq 0$, and $A_2 = 0$ otherwise.*

PROOF OF THEOREM 5. Without loss of generality, we assume $h_1 = X_1$ and $h_2 = X_2$. Since both $Q_{X_1}(X_1, \ldots, X_n)$ and $Q_{X_2}(X_1, \ldots, X_n)$ are saturated OAs, the GMA optimality and the overall $A_2 = s^n - 1$ follow from Lemmas 3, 5 and 6.

(i) The columns of $Q_{X_1}(X_1, \ldots, X_n)$ fall into three types: (a) $X_1$, (b) $X_1^2 + a_1X_1 + X_2$ and (c) $X_1^2 + a_1X_1 + b_1X_2 + g_1$, where $a_1, b_1 \in F_s$ and $g_1 \in H(X_3, \ldots, X_n)$. Similarly, the columns of $Q_{X_2}(X_1, \ldots, X_n)$ fall into three types: (a) $X_2$, (b) $X_2^2 + a_2X_2 + X_1$ and (c) $X_2^2 + a_2X_2 + b_2X_1 + g_2$, where $a_2, b_2 \in F_s$ and $g_2 \in H(X_3, \ldots, X_n)$. The projected $A_2$ values of all possible pairs can be found in Lemmas 11(ii), 11(i), 12(iv), 12(ii) and 12(i)(iii), respectively. In summary, we have the following aliasing patterns:

|  | $X_2$ | $X_2^2 + a_2X_2 + X_1$ | $X_2^2 + a_2X_2 + b_2X_1 + g_2$ |
| --- | --- | --- | --- |
| $X_1$ | 0 | $(s-1)/s$ | 0 |
| $X_1^2 + a_1X_1 + X_2$ | $(s-1)/s$ | $(s-1)^2/s^2$ | 0 |
| $X_1^2 + a_1X_1 + b_1X_2 + g_1$ | 0 | 0 | $\delta_{g_1,g_2}(s-1)/s^2$ |

where $\delta_{g_1,g_2}$ is equal to 1 if $g_1$ and $g_2$ are dependent and 0 otherwise. Each type (c) column in $Q_{X_1}(X_1, \ldots, X_n)$ is partially aliased with $s^2$ type (c) columns in $Q_{X_2}(X_1, \ldots, X_n)$. The result follows from the fact that the numbers of columns for each type are (a) 1, (b) $s$ and (c) $(s^n - s^2)/(s - 1)$, respectively.

(ii) From Lemmas 11 and 12, the possible projected $A_2$ values are 0, 1, 2 or 3.

(iii) From Lemmas 11 and 12, the possible projected $A_2$ values are 0, 1 or 3. Lemma 12(vi) shows that there is one fully aliased pair: $X_1^2 + X_2$ and $X_2^2 + X_1$, which has projected $A_2 = 3$. Since the overall $A_2 = 4^n - 1$, there must be $4^n - 4$ pairs with projected $A_2 = 1$. □

PROOF OF THEOREM 6. It follows from Theorem 5. □

PROOF OF THEOREM 7. We only need prove the GMA optimality. Since all linear functions form a saturated OA, the number of coincidences between any pair of rows of the resulting SSD is constant when $k = (s^n - 1)/(s - 1)$ and differs by at most one when $k = (s^n - 1)/(s - 1) - 1$. Therefore, the GMA optimality follows from Lemma 3. □

LEMMA 13. *Let $G \subset F_s$ and $|G| = k$. Suppose $X_1$ takes on values from $G$ only and all other $X_i$ take on values from $F_s$. Suppose $h_1, h_2 \in L(X_2, \ldots, X_n)$ and $a_1, a_2 \in F_s$.*

(i) *If $h_1$ and $h_2$ are independent, $a_1 X_1 + h_1$ and $a_2 X_1 + h_2$ are orthogonal.*

(ii) *If $h_1 = h_2$ and $a_1 \neq a_2$, the pair of columns $a_1 X_1 + h_1$ and $a_2 X_1 + h_2$ has projected $A_2 = (s - k)/k$.*

PROOF OF THEOREM 8. Without loss of generality, take $X_1$ as the branching column. The columns are $aX_1 + h$, where $a \in F_s$ and $h \in H(X_2, \ldots, X_n)$. By Lemma 13, each column is partially aliased with $s - 1$ columns with projected $A_2 = (s - k)/k$ and orthogonal to all other columns. Since there are $(s^n - s)/(s - 1)$ columns, there are $(s^n - s)/2$ pairs of nonorthogonal columns with projected $A_2 = (s - k)/k$. Therefore, the overall $A_2 = (s^n - s)(s - k)/(2k)$. Finally, the GMA optimality follows from Lemmas 3 and 5. □

LEMMA 14. *Let $G \subset F_s$ and $|G| = k$. Take $X_1^2 + X_2$ as the branching column of $Q_1(X_1, \ldots, X_n)$, that is, suppose all $X_i$, $i \neq 2$, take on values from $F_s$ and $X_1^2 + X_2$ takes on values from $G$ only. Suppose $h \in H(X_3, \ldots, X_n)$ and $a_1, a_2, b_1, b_2 \in F_s$.*

(i) *The pair of columns $X_1$ and $X_1^2 + a_1 X_1 + X_2$ has projected $A_2 = (s - k)/k$.*

(ii) *If $a_1 \neq a_2$, the pair of columns $X_1^2 + a_1 X_1 + X_2$ and $X_1^2 + a_2 X_1 + X_2$ has projected $A_2 = (s - k)/k$.*

(iii) *For $s$ odd, if $b_1 \neq b_2$, the pair of columns $X_1^2 + a_1 X_1 + b_1 X_2 + h$ and $X_1^2 + a_2 X_1 + b_2 X_2 + h$ has projected $A_2 = (s - k)/(ks)$.*

TABLE 1
*Some optimal three-level supersaturated designs*

| $N$ | $m$ | \multicolumn{5}{c|}{Projected $A_2$ values} | Source |
| --- | --- | --- | --- | --- | --- | --- | --- |
|  |  | 1/6 | 2/9 | 4/9 | 1/2 | 2/3 |  |
| 6 | 3 |  |  |  | 3 |  | Theorem 8, $n=2, k=2$ |
| 9 | 7 |  |  |  |  | 9 | Theorem 4, $n=2$ |
| 9 | 12 |  |  | 54 |  |  | Theorem 7, $n=2, k=4$ |
| 9 | 16 |  |  | 54 |  | 36 | Theorem 6, $n=2, k=4$ |
| 18 | 12 |  |  |  | 12 |  | Theorem 8, $n=3, k=2$ |
| 18 | 12 | 27 |  |  | 3 |  | Theorem 9, $n=3, k=2$ |
| 27 | 25 |  |  |  |  | 36 | Theorem 4, $n=3$ |
| 27 | 26 |  | 81 | 9 |  | 6 | Theorem 6, $n=3, k=2$ |
| 27 | 156 |  | 6318 | 702 |  |  | Theorem 7, $n=3, k=13$ |
| 27 | 169 |  | 6318 | 702 |  | 468 | Theorem 6, $n=3, k=13$ |
| 54 | 39 |  |  |  | 39 |  | Theorem 8, $n=4, k=2$ |
| 54 | 39 | 108 |  |  | 3 |  | Theorem 9, $n=4, k=2$ |

(iv) *For $s$ even, if $b_1 \neq b_2$ and $a_1 \neq a_2$, the pair of columns $X_1^2 + a_1 X_1 + b_1 X_2 + h$ and $X_1^2 + a_2 X_1 + b_2 X_2 + h$ has projected $A_2 \leq 1$.*

(v) *For $s = 4$, if $b_1 \neq b_2$ and $a_1 \neq a_2$, the pair of columns $X_1^2 + a_1 X_1 + b_1 X_2 + h$ and $X_1^2 + a_2 X_1 + b_2 X_2 + h$ has projected $A_2 = 0$ or $1$ for $k = 2$, and projected $A_2 = 1/9$ for $k = 3$.*

PROOF OF THEOREM 9. The GMA optimality follows from Lemmas 3 and 5. Since both designs in Theorems 8 and 9 have GMA, they must have the same overall $A_2 = (s^n - s)(s-k)/(2k)$.

(i) The columns of $Q_1(X_1, \ldots, X_n)$ are $X_1$, $X_1^2 + aX_1 + X_2$ and $X_1^2 + aX_1 + bX_2 + h$, where $a, b \in F_s$ and $h \in H(X_3, \ldots, X_n)$. By Lemma 14(i), the pair of columns $X_1$ and $X_1^2 + aX_1 + X_2$ has projected $A_2 = (s-k)/k$ when $a \neq 0$, and there are $s - 1$ such pairs; by Lemma 14(ii), the pair of columns $X_1^2 + a_1 X_1 + X_2$ and $X_1^2 + a_2 X_1 + X_2$ has projected $A_2 = (s-k)/k$ when $a_1 \neq a_2$, and there are $\binom{s-1}{2}$ such pairs since column $X_1^2 + X_2$ is removed; and by Lemma 14(iii), the pair of columns $X_1^2 + a_1 X_1 + b_1 X_2 + h$ and $X_1^2 + a_2 X_1 + b_2 X_2 + h$ has projected $A_2 = (s-k)/(ks)$ when $b_1 \neq b_2$, and there are $s^2 \binom{s}{2}(s^{n-2} - 1)/(s-1) = s(s^n - s^2)/2$ such pairs. It is easy to verify that all other pairs of columns are orthogonal.

(ii) and (iii) The proofs are similar to (i) and are omitted. □

**6. Some small designs and comparison.** Applying the construction methods, we can obtain many optimal multi-level SSDs. Tables 1–3 list the frequencies of nonzero projected $A_2$ values for some optimal 3-, 4- and 5-level

SSDs. All SSDs have the property that the number of coincidences between any pair of rows differs from each other by at most one; therefore, their overall $A_2$ values achieve the lower bound in Theorem 1 and they are optimal under GMA.

When $s = 4$ and $n = 2$, according to Theorem 6, the column juxtaposition of all five saturated OAs has 10 pairs of fully aliased columns. After removing one column from each pair, we obtain 15 columns with projected $A_2 = 0$ or 1. It can be verified that the resulting SSD has an overall $A_2$ value of 45 and achieves the lower bound in Theorem 1; therefore, it is optimal under GMA. Similarly, when $s = 4$ and $n = 3$, the column juxtaposition of all 21 saturated OAs has 210 pairs of fully aliased columns. After removing one column from each pair, we obtain 231 columns with projected $A_2 = 0$ or 1. It can be verified that the resulting SSD has an overall $A_2$ value of 3465 and achieves the lower bound in Theorem 1; therefore, it is also optimal under GMA.

We compare our SSDs based on Theorems 6 and 7 with existing designs from [2, 13, 22]. Since most designs are optimal under overall $A_2$ [or ave($\chi^2$)], we compare designs in terms of max($\chi^2$), ave($f$) and max($f$). Tables 4 and 5 show the comparisons for $N = 9, 16, 25$ and 27 runs in terms of ave($f$) and max($f$). For SSDs from Theorem 7, the first $m$ columns are used to evaluate these criteria. It is possible to find better designs if other columns are chosen.

Tables 4 and 5 indicate that our SSDs are competitive in terms of max($f$) but less competitive in terms of ave($f$). For $N = 9, 25$ and 27, SSDs based on Theorem 7 are better than existing ones in terms of both max($\chi^2$) and max($f$). In terms of ave($f$), our SSDs are worse than existing ones for $N = 9, 25$ but better for $N = 27$. For $N = 16$, our SSDs are less competitive.

Table 2
*Some optimal four-level supersaturated designs*

| | | Projected $A_2$ values | | | |
|---|---|---|---|---|---|
| $N$ | $m$ | 1/9 | 1/3 | 1 | Source |
| 8 | 4 | | | 6 | Theorem 8, $n = 2, k = 2$ |
| 12 | 4 | | 6 | | Theorem 8, $n = 2, k = 3$ |
| 16 | 9 | | | 12 | Theorem 4, $n = 2$ |
| 16 | 15 | | | 45 | Theorem $6^a$, $n = 2, k = 5$ |
| 32 | 20 | | | 30 | Theorem 8, $n = 3, k = 2$ |
| 48 | 20 | | 30 | | Theorem 8, $n = 3, k = 3$ |
| 48 | 20 | 72 | 6 | | Theorem 9, $n = 3, k = 3$ |
| 64 | 41 | | | 60 | Theorem 4, $n = 3$ |
| 64 | 231 | | | 3465 | Theorem $6^a$, $n = 3, k = 21$ |

[a] The design is obtained by removing fully aliased columns.

TABLE 3
*Some optimal five-level supersaturated designs*

| | | Projected $A_2$ values | | | | | | | |
|---|---|---|---|---|---|---|---|---|---|
| $N$ | $m$ | 2/15 | 1/4 | 3/10 | 16/25 | 2/3 | 4/5 | 3/2 | Source |
| 10 | 5 | | | | | | | 10 | Theorem 8, $n=2, k=2$ |
| 15 | 5 | | | | | 10 | | | Theorem 8, $n=2, k=3$ |
| 20 | 5 | | 10 | | | | | | Theorem 8, $n=2, k=4$ |
| 25 | 11 | | | | | | 25 | | Theorem 4, $n=2$ |
| 25 | 30 | | | | 375 | | | | Theorem 7, $n=2, k=6$ |
| 25 | 36 | | | | 375 | | 150 | | Theorem 6, $n=2, k=6$ |
| 50 | 30 | | | | | | | 60 | Theorem 8, $n=3, k=2$ |
| 50 | 30 | | | 250 | | | | 10 | Theorem 9, $n=3, k=2$ |
| 75 | 30 | | | | | 60 | | | Theorem 8, $n=3, k=3$ |
| 75 | 30 | 250 | | | | 10 | | | Theorem 9, $n=3, k=3$ |

TABLE 4
*Comparison of supersaturated designs in terms of* ave($f$)

| | | | Authors | | | | Aggarwal |
|---|---|---|---|---|---|---|---|
| $N$ | $s$ | $m$ | Theorem 6 | Theorem 7 | Fang et al. | Lu and Sun | and Gupta |
| 9 | 3 | 8 | 2.57 | 3.00 | 2.57 | 2.57 | 2.43 |
| 9 | 3 | 12 | 3.27 | 3.27 | 3.27 | 3.27 | 3.06 |
| 9 | 3 | 16 | 3.60 | | 3.60 | 3.60 | |
| 9 | 3 | 28 | | | | 4.00 | |
| 16 | 4 | 10 | 6.04 | | 4.36 | 4.84 | 4.93 |
| 16 | 4 | 15 | 6.86 | | 5.60 | 6.23 | 6.27 |
| 16 | 4 | 20 | | | 6.25 | 6.95 | |
| 16 | 4 | 40 | | | | 7.87 | |
| 25 | 5 | 12 | 8.33 | 9.55 | 6.42 | 8.06 | 7.45 |
| 25 | 5 | 18 | 10.78 | 10.98 | 8.41 | 10.42 | 9.52 |
| 25 | 5 | 24 | 11.96 | 11.67 | 10.20 | 11.66 | 10.86 |
| 25 | 5 | 30 | 12.64 | 12.07 | | 12.33 | 11.33 |
| 25 | 5 | 36 | 13.10 | | | 12.73 | |
| 27 | 3 | 26 | 3.66 | 3.77 | 3.78 | | 4.26 |
| 27 | 3 | 39 | 4.81 | 4.81 | 5.27 | | 5.63 |
| 27 | 3 | 52 | 5.38 | 5.97 | 5.98 | | 6.32 |
| 27 | 3 | 65 | 5.71 | 6.28 | | | 6.73 |
| 27 | 3 | 156 | 6.49 | 6.97 | | | |
| 27 | 3 | 169 | 6.53 | | | | |

We also compare our SSDs based on Theorems 8 and 9 with designs from Lu, Hu and Zheng [21], who constructed some small SSDs based on resolvable balanced incomplete block designs. We find that our designs have the same $\max(\chi^2)$ and $\max(f)$ values as theirs. Note that their methods depend on the existence of resolvable balanced incomplete block designs, which is not an easy task itself. In contrast, our algebraic construction is general and works for any $s$ and $n$. Indeed, SSDs based on Theorems 8 and 9 with $n \geq 3$ are not available in [21].

TABLE 5
*Comparison of supersaturated designs in terms of* $\max(f)$

| $N$ | $s$ | $m$ | Authors | | | | Aggarwal |
| --- | --- | --- | --- | --- | --- | --- | --- |
| | | | Theorem 6 | Theorem 7 | Fang et al. | Lu and Sun | and Gupta |
| 9 | 3 | 8 | 6 | 4 | 6 | 6 | 8 |
| 9 | 3 | 12 | 6 | 4 | 6 | 6 | 8 |
| 9 | 3 | 16 | 6 | | 6 | 6 | |
| 9 | 3 | 28 | | | | 6 | |
| 16 | 4 | 10 | 16 | | 12 | 12 | 12 |
| 16 | 4 | 15 | 16 | | 12 | 12 | 14 |
| 16 | 4 | 20 | | | 16 | 12 | |
| 16 | 4 | 40 | | | | 16 | |
| 25 | 5 | 12 | 20 | 14 | 22 | 18 | 24 |
| 25 | 5 | 18 | 20 | 14 | 24 | 20 | 24 |
| 25 | 5 | 24 | 20 | 14 | 30 | 20 | 32 |
| 25 | 5 | 30 | 20 | 14 | | 22 | 32 |
| 25 | 5 | 36 | 20 | | | 22 | |
| 27 | 3 | 26 | 18 | 12 | 16 | | 16 |
| 27 | 3 | 39 | 18 | 12 | 18 | | 16 |
| 27 | 3 | 52 | 18 | 12 | 18 | | 16 |
| 27 | 3 | 65 | 18 | 12 | | | 16 |
| 27 | 3 | 156 | 18 | 12 | | | |
| 27 | 3 | 169 | 18 | | | | |

We have presented several classes of optimal SSDs whose columns are represented by linear and quadratic polynomials and analytically studied the aliasing structure among columns. SSDs based on Theorem 7 are generally preferred to those based on Theorems 4 and 6, because the former have smaller maximum pairwise aliasing in terms of both $\max(\chi^2)$ and $\max(f)$. Nevertheless, SSDs based on Theorem 4 are useful in some situations. For example, when the experimenter isolates one important factor and wants to estimate it efficiently, then that factor should be assigned to column $X_1$ which is orthogonal to all other columns.

For easy reference, the Appendix lists some small SSDs based on Theorems 6 and 7. All SSDs in Tables 1–3 are explicitly listed online at www.stat.ucla.edu/~hqxu/.

**7. Mixed-level SSDs.** The GMA criterion works for mixed-level SSDs. The following lemma shows that $\mathrm{ave}(\chi^2)$ is again equivalent to $A_2$ for mixed-level SSDs. The proof is similar to Lemma 4(i) and is thus omitted.

LEMMA 15. *For an $SSD(N, s_1 s_2 \cdots s_m)$, $\mathrm{ave}(\chi^2) = NA_2/[m(m-1)/2]$.*

The following lower bound of $A_2$ for mixed-level SSDs generalizes Lemma 2.

THEOREM 10. *For an $SSD(N, s_1 s_2 \cdots s_m)$, $A_2 \geq (\sum_{k=1}^m s_k - m)(\sum_{k=1}^m s_k - m - N + 1)/[2(N-1)]$.*

PROOF. Let $D = [x_{ik}]$ be the $N \times m$ design matrix. Let $\delta_{ij}^{(w)}(D) = \sum_{k=1}^m s_k \delta(x_{ik}, x_{jk})$ be the weighted number of coincidences between the rows $i$ and $j$, where $\delta(x, y) = 1$ if $x = y$ and 0 otherwise. Define $J_t(D) = \sum_{1 \leq i < j \leq N} [\delta_{ij}^{(w)}(D)]^t$ for $t = 1, 2$. Xu ([32], Lemma 2) showed that $J_2(D) = N^2 A_2(D) + N[Nm(m-1) + N \sum s_k - (\sum s_k)^2]/2$. The proof there also implies that $J_1(D) = N(Nm - \sum s_k)/2$. Thus, the lower bound of $A_2(D)$ follows from the Cauchy–Schwarz inequality $[N(N-1)/2]J_2(D) \geq [J_1(D)]^2$ and some straightforward algebra. $\square$

The lower bound in Theorem 10 is equivalent to the lower bound of $\mathrm{ave}(\chi^2)$ given by Yamada and Matsui ([38], Theorem 1) and the lower bound of $K_2$ given by Xu ([33], Theorem 6).

The situation for mixed-level SSDs is more complicated than that for multi-level SSDs. Further development is needed to learn whether mixed-level SSDs achieving the lower bound in Theorem 10 are optimal under GMA.

Finally, optimal mixed-level SSDs can be generated from multi-level SSDs via the method of replacement. For the method of replacement, see [15] and Wu and Hamada ([31], Section 7.7). Xu and Wu [34] showed that the generalized wordlength patterns are invariant with respect to the choice of orthonormal contrasts. As a result, when one or more $s$-level columns are replaced with a saturated OA of run size $s$ and strength 2, the resulting mixed-level SSD has the same overall $A_2$ values as the original multi-level SSD, and the maximum projected $A_2$ value of the mixed-level SSD is always less than or equal to that of the original multi-level SSD. Furthermore, if the original multi-level SSD achieves the lower bound of $A_2$ in Lemma 2, then the mixed-level SSD achieves the lower bound of $A_2$ in Theorem 10 and thus is $A_2$ optimal.

The following example illustrates these ideas. The $SSD(81, 9^{100})$ from Theorem 6 (with $n = 2$, $s = 9$ and $k = 10$) has overall $A_2 = 3600$ and maximum projected $A_2 = 8/9$. We obtain an $SSD(81, 9^{100-i} 3^{4i})$ by replacing $i$ 9-level columns with four 3-level columns that form an $OA(9, 4, 3, 2)$ for $1 \leq i < 100$. All these mixed-level SSDs have overall $A_2 = 3600$ and projected $A_2 \leq 8/9$. In addition, all these mixed-level SSDs achieve the lower bound of $A_2$ in Theorem 10; thus, they are optimal under both $A_2$ and $\mathrm{ave}(\chi^2)$. However, we do not know whether they are optimal under GMA.

# APPENDIX

TABLE 6
$SSD(9,3^{16})$ and $SSD(9,3^{12})$ via Theorems 6 and 7 ($n=2$, $s=3$, $k=4$)

| Run | 1 | 2 | 3 | 4 | 5 | 6 | 7 | 8 | 9 | 10 | 11 | 12 | 13 | 14 | 15 | 16 |
|---|---|---|---|---|---|---|---|---|---|---|---|---|---|---|---|---|
| 1 | 0 | 0 | 0 | 0 | 0 | 0 | 0 | 0 | 0 | 0 | 0 | 0 | 0 | 0 | 0 | 0 |
| 2 | 0 | 1 | 1 | 1 | 1 | 1 | 2 | 0 | 1 | 1 | 2 | 0 | 1 | 1 | 2 | 0 |
| 3 | 0 | 2 | 2 | 2 | 2 | 1 | 0 | 2 | 2 | 1 | 0 | 2 | 2 | 1 | 0 | 2 |
| 4 | 1 | 1 | 2 | 0 | 0 | 1 | 1 | 1 | 1 | 2 | 0 | 1 | 2 | 2 | 1 | 0 |
| 5 | 1 | 2 | 0 | 1 | 1 | 2 | 0 | 1 | 2 | 2 | 1 | 0 | 0 | 1 | 1 | 1 |
| 6 | 1 | 0 | 1 | 2 | 2 | 2 | 1 | 0 | 0 | 1 | 1 | 1 | 1 | 2 | 0 | 1 |
| 7 | 2 | 1 | 0 | 2 | 0 | 2 | 2 | 2 | 2 | 0 | 2 | 1 | 1 | 0 | 1 | 2 |
| 8 | 2 | 2 | 1 | 0 | 1 | 0 | 1 | 2 | 0 | 2 | 2 | 2 | 2 | 0 | 2 | 1 |
| 9 | 2 | 0 | 2 | 1 | 2 | 0 | 2 | 1 | 1 | 0 | 1 | 2 | 0 | 2 | 2 | 2 |

NOTES. (i) All columns form an $SSD(9,3^{16})$ with overall $A_2 = 48$; the possible projected $A_2$ values are 0, 4/9 and 2/3 with frequencies 30, 54 and 36. (ii) Removing columns (1, 5, 9, 13) yields an $SSD(9,3^{12})$ with overall $A_2 = 24$; the possible projected $A_2$ values are 0 and 4/9 with frequencies 12 and 54.

TABLE 7
$SSD(16,4^{15})$ via Theorem 6 ($n=2, s=4, k=5$)

| Run | 1 | 2 | 3 | 4 | 5 | 6 | 7 | 8 | 9 | 10 | 11 | 12 | 13 | 14 | 15 |
|---|---|---|---|---|---|---|---|---|---|---|---|---|---|---|---|
| 1 | 0 | 0 | 0 | 0 | 0 | 0 | 0 | 0 | 0 | 0 | 0 | 0 | 0 | 0 | 0 |
| 2 | 0 | 1 | 1 | 1 | 1 | 1 | 0 | 3 | 2 | 1 | 3 | 2 | 1 | 0 | 1 |
| 3 | 0 | 2 | 2 | 2 | 2 | 2 | 1 | 0 | 2 | 2 | 0 | 2 | 2 | 1 | 2 |
| 4 | 0 | 3 | 3 | 3 | 3 | 3 | 1 | 3 | 0 | 3 | 3 | 0 | 3 | 1 | 3 |
| 5 | 1 | 1 | 0 | 3 | 2 | 0 | 1 | 1 | 1 | 1 | 2 | 3 | 2 | 0 | 3 |
| 6 | 1 | 0 | 1 | 2 | 3 | 1 | 1 | 2 | 3 | 0 | 1 | 1 | 3 | 0 | 2 |
| 7 | 1 | 3 | 2 | 1 | 0 | 2 | 0 | 1 | 3 | 3 | 2 | 1 | 0 | 1 | 1 |
| 8 | 1 | 2 | 3 | 0 | 1 | 3 | 0 | 2 | 1 | 2 | 1 | 3 | 1 | 1 | 0 |
| 9 | 2 | 3 | 1 | 0 | 2 | 0 | 2 | 2 | 2 | 2 | 2 | 0 | 3 | 3 | 1 |
| 10 | 2 | 2 | 0 | 1 | 3 | 1 | 2 | 1 | 0 | 3 | 1 | 2 | 2 | 3 | 0 |
| 11 | 2 | 1 | 3 | 2 | 0 | 2 | 3 | 2 | 0 | 0 | 2 | 2 | 1 | 2 | 3 |
| 12 | 2 | 0 | 2 | 3 | 1 | 3 | 3 | 1 | 2 | 1 | 1 | 0 | 0 | 2 | 2 |
| 13 | 3 | 2 | 1 | 3 | 0 | 0 | 3 | 3 | 3 | 3 | 0 | 3 | 1 | 3 | 2 |
| 14 | 3 | 3 | 0 | 2 | 1 | 1 | 3 | 0 | 1 | 2 | 3 | 1 | 0 | 3 | 3 |
| 15 | 3 | 0 | 3 | 1 | 2 | 2 | 2 | 3 | 1 | 1 | 0 | 1 | 3 | 2 | 0 |
| 16 | 3 | 1 | 2 | 0 | 3 | 3 | 2 | 0 | 3 | 0 | 3 | 3 | 2 | 2 | 1 |

NOTES. (i) This design is obtained by removing a column from each pair of fully aliased columns. (ii) The overall $A_2$ value is 45; the possible projected $A_2$ values are 0 and 1 with frequencies 60 and 45.

TABLE 8
$SSD(25, 5^{36})$ and $SSD(25, 5^{30})$ via Theorems 6 and 7 ($n=2, s=5, k=6$)

| Run | 1 | 2 | 3 | 4 | 5 | 6 | 7 | 8 | 9 | 10 | 11 | 12 | 13 | 14 | 15 | 16 | 17 | 18 | 19 | 20 | 21 | 22 | 23 | 24 | 25 | 26 | 27 | 28 | 29 | 30 | 31 | 32 | 33 | 34 | 35 | 36 |
|---|---|---|---|---|---|---|---|---|---|---|---|---|---|---|---|---|---|---|---|---|---|---|---|---|---|---|---|---|---|---|---|---|---|---|---|---|
| 1 | 0 | 0 | 0 | 0 | 0 | 0 | 0 | 0 | 0 | 0 | 0 | 0 | 0 | 0 | 0 | 0 | 0 | 0 | 0 | 0 | 0 | 0 | 0 | 0 | 0 | 0 | 0 | 0 | 0 | 0 | 0 | 0 | 0 | 0 | 0 | 0 |
| 2 | 0 | 1 | 1 | 1 | 1 | 1 | 1 | 1 | 2 | 3 | 4 | 0 | 1 | 1 | 2 | 3 | 4 | 0 | 1 | 1 | 2 | 3 | 4 | 0 | 1 | 1 | 2 | 3 | 4 | 0 | 1 | 1 | 2 | 3 | 4 | 0 |
| 3 | 0 | 2 | 2 | 2 | 2 | 2 | 2 | 4 | 1 | 3 | 0 | 2 | 2 | 4 | 1 | 3 | 0 | 2 | 2 | 4 | 1 | 3 | 0 | 2 | 2 | 4 | 1 | 3 | 0 | 2 | 2 | 4 | 1 | 3 | 0 | 2 |
| 4 | 0 | 3 | 3 | 3 | 3 | 3 | 3 | 4 | 2 | 0 | 3 | 1 | 3 | 4 | 2 | 0 | 3 | 1 | 3 | 4 | 2 | 0 | 3 | 1 | 3 | 4 | 2 | 0 | 3 | 1 | 3 | 4 | 2 | 0 | 3 | 1 |
| 5 | 0 | 4 | 4 | 4 | 4 | 4 | 4 | 1 | 0 | 4 | 3 | 2 | 4 | 1 | 0 | 4 | 3 | 2 | 4 | 1 | 0 | 4 | 3 | 2 | 4 | 1 | 0 | 4 | 3 | 2 | 4 | 1 | 0 | 4 | 3 | 2 |
| 6 | 1 | 1 | 2 | 3 | 4 | 0 | 0 | 1 | 1 | 1 | 1 | 1 | 2 | 3 | 4 | 0 | 1 | 2 | 0 | 2 | 4 | 1 | 3 | 3 | 0 | 3 | 1 | 4 | 2 | 4 | 2 | 1 | 0 | 4 | 3 | 3 |
| 7 | 1 | 2 | 3 | 4 | 0 | 1 | 1 | 2 | 3 | 4 | 0 | 1 | 2 | 0 | 2 | 4 | 1 | 3 | 3 | 0 | 3 | 1 | 4 | 2 | 4 | 2 | 1 | 0 | 4 | 3 | 0 | 1 | 1 | 1 | 1 | 1 |
| 8 | 1 | 3 | 4 | 0 | 1 | 2 | 2 | 0 | 2 | 4 | 1 | 3 | 3 | 0 | 3 | 1 | 4 | 2 | 4 | 2 | 1 | 0 | 4 | 3 | 0 | 1 | 1 | 1 | 1 | 1 | 2 | 3 | 4 | 0 | 1 |
| 9 | 1 | 4 | 0 | 1 | 2 | 3 | 3 | 0 | 3 | 1 | 4 | 2 | 4 | 2 | 1 | 0 | 4 | 3 | 0 | 1 | 1 | 1 | 1 | 1 | 2 | 3 | 4 | 0 | 1 | 2 | 0 | 2 | 4 | 1 | 3 |
| 10 | 1 | 0 | 1 | 2 | 3 | 4 | 4 | 2 | 1 | 0 | 4 | 3 | 0 | 1 | 1 | 1 | 1 | 1 | 2 | 3 | 4 | 0 | 1 | 2 | 0 | 2 | 4 | 1 | 3 | 3 | 0 | 3 | 1 | 4 | 2 |
| 11 | 2 | 4 | 1 | 3 | 0 | 2 | 0 | 2 | 2 | 2 | 2 | 2 | 1 | 3 | 0 | 2 | 4 | 4 | 3 | 2 | 1 | 0 | 4 | 1 | 3 | 4 | 0 | 1 | 2 | 3 | 1 | 4 | 2 | 0 | 3 |
| 12 | 2 | 0 | 2 | 4 | 1 | 3 | 1 | 3 | 4 | 0 | 1 | 2 | 3 | 1 | 4 | 2 | 0 | 3 | 0 | 2 | 2 | 2 | 2 | 2 | 1 | 3 | 0 | 2 | 4 | 4 | 3 | 2 | 1 | 0 | 4 |
| 13 | 2 | 1 | 3 | 0 | 2 | 4 | 2 | 1 | 3 | 0 | 2 | 4 | 4 | 3 | 2 | 1 | 0 | 4 | 1 | 3 | 4 | 0 | 1 | 2 | 3 | 1 | 4 | 2 | 0 | 3 | 0 | 2 | 2 | 2 | 2 | 2 |
| 14 | 2 | 2 | 4 | 1 | 3 | 0 | 3 | 1 | 4 | 2 | 0 | 3 | 0 | 2 | 2 | 2 | 2 | 2 | 1 | 3 | 0 | 2 | 4 | 4 | 3 | 2 | 1 | 0 | 4 | 1 | 3 | 4 | 0 | 1 | 2 |
| 15 | 2 | 3 | 0 | 2 | 4 | 1 | 4 | 3 | 2 | 1 | 0 | 4 | 1 | 3 | 4 | 0 | 1 | 2 | 3 | 1 | 4 | 2 | 0 | 3 | 0 | 2 | 2 | 2 | 2 | 2 | 1 | 3 | 0 | 2 | 4 |
| 16 | 3 | 4 | 2 | 0 | 3 | 1 | 0 | 3 | 3 | 3 | 3 | 3 | 2 | 0 | 3 | 1 | 4 | 1 | 4 | 0 | 1 | 2 | 3 | 4 | 4 | 3 | 2 | 1 | 0 | 2 | 2 | 4 | 1 | 3 | 0 |
| 17 | 3 | 0 | 3 | 1 | 4 | 2 | 1 | 4 | 0 | 1 | 2 | 3 | 4 | 4 | 3 | 2 | 1 | 0 | 2 | 2 | 4 | 1 | 3 | 0 | 0 | 3 | 3 | 3 | 3 | 3 | 2 | 0 | 3 | 1 | 4 |
| 18 | 3 | 1 | 4 | 2 | 0 | 3 | 2 | 2 | 4 | 1 | 3 | 0 | 0 | 3 | 3 | 3 | 3 | 3 | 2 | 0 | 3 | 1 | 4 | 1 | 4 | 0 | 1 | 2 | 3 | 4 | 4 | 3 | 2 | 1 | 0 |
| 19 | 3 | 2 | 0 | 3 | 1 | 4 | 3 | 2 | 0 | 3 | 1 | 4 | 0 | 1 | 2 | 3 | 4 | 4 | 3 | 2 | 1 | 0 | 2 | 2 | 4 | 1 | 3 | 0 | 0 | 3 | 3 | 3 | 3 | 3 | 3 |
| 20 | 3 | 3 | 1 | 4 | 2 | 0 | 4 | 4 | 3 | 2 | 1 | 0 | 2 | 2 | 4 | 1 | 3 | 0 | 0 | 3 | 3 | 3 | 3 | 3 | 2 | 0 | 3 | 1 | 4 | 1 | 4 | 0 | 1 | 2 | 3 |
| 21 | 4 | 1 | 0 | 4 | 3 | 2 | 0 | 4 | 4 | 4 | 4 | 4 | 4 | 0 | 4 | 3 | 2 | 1 | 3 | 3 | 1 | 4 | 2 | 0 | 2 | 3 | 0 | 2 | 4 | 1 | 1 | 0 | 1 | 2 | 3 | 4 |
| 22 | 4 | 2 | 1 | 0 | 4 | 3 | 1 | 0 | 1 | 2 | 3 | 4 | 0 | 4 | 4 | 4 | 4 | 4 | 4 | 0 | 4 | 3 | 2 | 1 | 3 | 3 | 1 | 4 | 2 | 0 | 2 | 3 | 0 | 2 | 4 | 1 |
| 23 | 4 | 3 | 2 | 1 | 0 | 4 | 2 | 3 | 0 | 2 | 4 | 1 | 1 | 0 | 1 | 2 | 3 | 4 | 0 | 4 | 4 | 4 | 4 | 4 | 0 | 4 | 3 | 2 | 1 | 3 | 3 | 1 | 4 | 2 | 0 |
| 24 | 4 | 4 | 3 | 2 | 1 | 0 | 3 | 3 | 1 | 4 | 2 | 0 | 2 | 3 | 0 | 2 | 4 | 1 | 1 | 0 | 1 | 2 | 3 | 4 | 0 | 4 | 4 | 4 | 4 | 4 | 0 | 4 | 3 | 2 | 1 |
| 25 | 4 | 0 | 4 | 3 | 2 | 1 | 4 | 0 | 4 | 3 | 2 | 1 | 3 | 3 | 1 | 4 | 2 | 0 | 2 | 3 | 0 | 2 | 4 | 1 | 1 | 0 | 1 | 2 | 3 | 4 | 0 | 4 | 4 | 4 | 4 | 4 |

NOTES. (i) All columns form an $SSD(25, 5^{36})$ with overall $A_2 = 360$; the possible projected $A_2$ values are 0, 16/25 and 4/5 with frequencies 105, 375 and 150. (ii) Removing columns (1, 7, 13, 19, 25, 31) yields an $SSD(25, 5^{30})$ with overall $A_2 = 240$; the possible projected $A_2$ values are 0 and 16/25 with frequencies 60 and 375.

**Acknowledgments.** The authors thank an Associate Editor and the referees for their valuable comments.

| | |
|---|---|
| Department of Statistics | School of Industrial |
| University of California | and Systems Engineering |
| Los Angeles, California 90095-1554 | Georgia Institute of Technology |
| USA | Atlanta, Georgia 30332-0205 |
| E-mail: hqxu@stat.ucla.edu | USA |
| | E-mail: jeffwu@isye.gatech.edu |